\documentclass[titlepage,draft,12pt]{article} 
\usepackage{amsfonts,latexsym,pstricks,pst-plot,graphicx} 
\usepackage[a4paper]{geometry}
\newcommand{\allarga}{\rule[-3ex]{0ex}{6.5ex}}
\newcommand{\ep}{\varepsilon}

\newcommand{\cep}{c_{\ep}}
\newcommand{\rep}{r_{\ep}}
\newcommand{\roep}{\rho_{\ep}}
\newcommand{\uep}{u_{\ep}}
\newcommand{\Gep}{G_{\ep}}
\newcommand{\Pep}{P_{\ep}}
\newcommand{\Qep}{Q_{\ep}}
\newcommand{\tetep}{\theta_{\ep}}
\newcommand{\Hep}{\mathcal{H}_{\ep}}
\newcommand{\m}[1]{m(\au{#1}^{2})}
\newcommand{\au}[1]{|A^{1/2}#1|}
\newcommand{\auq}[1]{|A^{1/2}#1|^{2}}
\newcommand{\dat}{D(A^{3/2})}
\newcommand{\da}{D(A)}
\newcommand{\dau}{D(A^{1/2})}

\newtheorem{thm}{Theorem}[section]

\newtheorem{open}{Open problem}

\title{Hyperbolic-parabolic singular perturbation for Kirchhoff 
equations with weak dissipation}
\author{Marina Ghisi\medskip\\ {\normalsize
Universit\`a degli Studi di Pisa} \\{\normalsize Dipartimento di
Matematica ``Leonida Tonelli''}\\
{\normalsize 
PISA (Italy)}\\  
{\normalsize e-mail: \texttt{ghisi@dm.unipi.it}}\and 
Massimo Gobbino\medskip\\ {\normalsize Universit\`a degli Studi di Pisa} 
\\{\normalsize Dipartimento di Matematica Applicata ``Ulisse Dini''}\\ 
{\normalsize 
 PISA (Italy)}\\  
{\normalsize e-mail: \texttt{m.gobbino@dma.unipi.it}}}
\date{}

\begin{document}
\maketitle
\begin{abstract}
	We consider Kirchhoff equations with a small parameter $\ep$ in
	front of the second-order time-derivative, and a dissipative term
	whose coefficient may tend to $0$ as $t\to +\infty$ (weak
	dissipation).
	
	In this note we present some recent results concerning existence
	of global solutions, and their asymptotic behavior both as $t\to
	+\infty$ and as $\ep\to 0^{+}$.  Since the limit equation is of
	parabolic type, this is usually referred to as a
	hyperbolic-parabolic singular perturbation problem.
	
	We show in particular that the equation exhibits hyperbolic or
	parabolic behavior depending on the values of the parameters.
	
\vspace{1cm}

\noindent{\bf Mathematics Subject Classification 2000 (MSC2000):}
35B25, 35L70, 35B40.

\vspace{1cm} 

\end{abstract}
 
\section{Introduction}

Let $H$ be a separable real Hilbert space.  For every $x$ and $y$ in
$H$, $|x|$ denotes the norm of $x$, and $\langle x,y\rangle$ denotes
the scalar product of $x$ and $y$.  Let $A$ be a self-adjoint linear
operator on $H$ with dense domain $D(A)$.  We assume that $A$ is
nonnegative, namely $\langle Ax,x\rangle\geq 0$ for every $x\in D(A)$,
so that for every $\alpha\geq 0$ the power $A^{\alpha}x$ is defined
provided that $x$ lies in a suitable domain $D(A^{\alpha})$.

Let $b:[0,+\infty)\to(0,+\infty)$ and $m:[0,+\infty)\to[0,+\infty)$ be
two given functions. For every $\ep>0$ we consider the 
Cauchy problem
\begin{equation}
	\ep\uep''(t)+b(t)\uep'(t)+
	\m{\uep(t)}A\uep(t)=0,
	\label{pbm:h-eq-gen}
\end{equation}
\begin{equation}
	u(0)=u_0,\hspace{3em}u'(0)=u_1.
	\label{pbm:h-data}
\end{equation}

This is the dissipative version of the celebrated equation introduced 
by G.\ Kirchhoff in~\cite{kirchhoff} as a simplified model for 
transversal vibrations of elastic strings. We refer to the 
survey~\cite{gg:survey-nd} for the non-dissipative case where 
$\ep=1$ and $b(t)\equiv 0$. Let us set
$$\mu:=\inf_{\sigma\geq 0}m(\sigma), 
\quad\quad 
\delta:=\inf_{t\geq 0}b(t),
\quad\quad 
\nu:=\inf\left\{\frac{\langle Ax,x\rangle}{|x|^{2}}:x\in D(A),\ x\neq
0\right\}.$$

Several features of (\ref{pbm:h-eq-gen}) depend on the values of 
$\mu$, $\delta$, $\nu$. Let us recall some standard terminology.
\begin{itemize}
	\item \emph{Non-degenerate vs degenerate equations}\quad These
	terms refer to the nonlinearity.  Equation (\ref{pbm:h-eq-gen}) is
	called \emph{nondegenerate} or \emph{strictly hyperbolic} if
	$\mu>0$, and \emph{degenerate} or \emph{weakly hyperbolic} if
	$\mu\geq 0$.  The Cauchy problem (\ref{pbm:h-eq-gen}),
	(\ref{pbm:h-data}) is called \emph{mildly degenerate} if $\mu\geq
	0$ but
	\begin{equation}
		\m{u_{0}}\neq 0.
		\label{hp:mdg}
	\end{equation}
	Whenever we consider degenerate equations, we always limit
	ourselves to the mildly degenerate case.  The \emph{really
	degenerate} case where $\m{u_{0}}=0$ seems to be still quite
	unexplored.

	\item \emph{Constant vs weak dissipation}\quad We have
	\emph{constant dissipation} when $b(t)=\delta>0$ for every $t\geq
	0$, and \emph{weak dissipation} when $b(t)\to 0$ as $t\to
	+\infty$.  Almost all known results for the constant dissipation
	case can be easily extended to non-constant dissipation
	coefficients provided that $\delta>0$ and $b'(t)$ is bounded.  For
	simplicity we often limit ourselves to the model case where
	$b(t)=(1+t)^{-p}$ for some $p\geq 0$, the case $p=0$ corresponding
	to constant dissipation.
	
	In this note we don't consider equations with \emph{strong
	dissipation}, which usually refers to dissipative terms of the
	form $A^{\alpha}\uep'(t)$ with $\alpha>0$, or better
	$\alpha\geq 1/2$.

	\item \emph{Coercive vs non-coercive operators}\quad The operator
	$A$ is called \emph{coercive} when $\nu>0$, and it is called
	\emph{noncoercive} when $\nu\geq 0$.  This property of the
	operator has a great influence on the asymptotic behavior of
	solutions.
\end{itemize}

The singular perturbation problem in its generality consists in
proving the convergence of solutions of (\ref{pbm:h-eq-gen}),
(\ref{pbm:h-data}) to solutions of the first order problem
\begin{equation}
	b(t)\uep'(t)+ \m{\uep(t)}A\uep(t)=0, \hspace{2em}
	u(0)=u_{0},
	\label{pbm:par}
\end{equation}
obtained setting formally $\ep=0$ in (\ref{pbm:h-eq-gen}), and
omitting the second initial condition in~(\ref{pbm:h-data}).
Following the approach introduced by \textsc{J.\ L.\
Lions}~\cite{lions} in the linear case, one defines the corrector
$\tetep(t)$ as the solution of the second order \emph{linear} problem
$$\ep\tetep''(t)+b(t)\tetep'(t)=0 \hspace{2em}
	\forall t\geq 0,$$
\begin{equation}
	\tetep(0)=0,\hspace{2em}\tetep'(0)=u_1+\frac{1}{b(0)}
	\m{u_{0}}Au_{0}=:w_{0}.
	\label{pbm:theta-data}
\end{equation}

It is easy to see that $\tetep'(0)=\uep'(0)-u'(0)$, hence this
corrector keeps into account the boundary layer due to the loss of one
initial condition.  Finally one defines $\rep(t)$ and $\roep(t)$ in
such a way that
\begin{equation}
	\uep(t)=u(t)+\tetep(t)+\rep(t)=u(t)+\roep(t)\quad\quad\forall
	t\geq 0.
	\label{defn:rep}
\end{equation}

With these notations, the singular perturbation problem consists in
proving that $\rep(t)\to 0$ or $\roep(t)\to 0$ in some sense as
$\ep\to 0^{+}$. The general problem can be split into at least six 
subproblems.
\begin{enumerate}
	\renewcommand{\labelenumi}{(\arabic{enumi})} 
	\item \emph{Parabolic problem: global existence and decay
	estimates}\quad This is the first and usually easiest step in the
	theory.  It consists in proving that (\ref{pbm:par}) admits a
	unique global solution $u(t)$, and then in estimating its decay
	rate as $t\to+\infty$.  This decay rate is afterwards used as a
	benchmark for the decay rate of solutions of the hyperbolic
	problem.

	\item \emph{Local existence for the hyperbolic problem and
	local-in-time error estimates}\quad Let $T>0$ be fixed.  This
	subproblem consists in proving that, for every $\ep>0$ small
	enough, the solution $\uep(t)$ of the hyperbolic problem
	(\ref{pbm:h-eq-gen}), (\ref{pbm:h-data}) is defined (at least) on
	the interval $[0,T]$, and in this interval $\uep(t)$ converges to
	the solution $u(t)$ of the limit problem. In this case the 
	smallness of $\ep$, as well as the convergence rates, may depend 
	on $T$.

	\item \emph{Hyperbolic problem: global existence}\quad This
	subproblem consists in proving that problem (\ref{pbm:h-eq-gen}),
	(\ref{pbm:h-data}) admits a global-in-time solution provided that
	$\ep>0$ is small enough.  From the point of view of existence, this
	is a strengthening of the previous step, and of course in general
	it requires stronger assumptions.
	
	Existence of global solutions without the smallness assumption on
	$\ep$ is a widely open question, which seems to be as difficult as
	the non-dissipative case (see Section~\ref{sec:hyperbolic}).

	\item \emph{Hyperbolic problem: decay estimates}\quad Once we know
	that the hyperbolic problem admits a global-in-time solution
	$\uep(t)$, a natural question concerns its behavior as $t\to
	+\infty$ ($\ep$ is now small and fixed).  What one expects in
	reasonable situations is that $\uep(t)$ decays as the solution
	$u(t)$ of the corresponding parabolic equation.  This is what has
	been actually proved in many cases.

	\item \emph{Singular perturbation problem: global-in-time error
	estimates}\quad This subproblem is just the global-in-time version
	of subproblem (2).  The goal is therefore to give
	\emph{time-independent} estimates on $\roep(t)$ or $\rep(t)$ as
	$\ep\to 0^{+}$.

	\item \emph{Singular perturbation problem: decay-error
	estimates}\quad This is the meeting point of subproblems (4) and
	(5), and it is the ultimate goal of the theory.  It consists in
	estimating in the same time the behavior of $\uep(t)$ as
	$t\to+\infty$ and as $\ep\to 0^{+}$.  
	The general form of a decay-error estimate is something like
	$$|A^{\alpha}\roep(t)|\leq\omega(\ep)\gamma(t) \quad\mbox{or}\quad
	|A^{\alpha}\rep(t)|\leq\omega(\ep)\gamma(t).$$
	
	Of course one expects $\gamma(t)$ to be the decay rate of
	solutions of the parabolic problem (or even better), and
	$\omega(\ep)$ to be the convergence rate which appears in the
	local-in-time error estimates.
\end{enumerate}

This program has generated a considerable literature in the last
thirty years, for which we refer to the introductions of the following
sections.  In this note we sum up the state of the art and the main
open questions.  A rough overview is provided by
Table~\ref{tab:overview}, where we show, under different assumptions,
which subproblems have received a reasonable or partial answer up to
now.  We focus in particular on the model dissipation coefficient of
the form $b(t):=(1+t)^{-p}$ with $p\geq 0$, and on nonlinear terms
which are either non-degenerate or of the form
$m(\sigma)=\sigma^{\gamma}$ for some $\gamma>0$ (note that we allow
also the non-Lipschitz case $\gamma\in(0,1))$.
\begin{table}[htp]
	\centering
	\renewcommand{\arraystretch}{1.4}
	\begin{tabular}{|c||c|c|c|}
		\hline
		& $p=0$ & $p\in(0,1]$, $\nu>0$ & $p\in(0,1]$, $\nu\geq 0$  \\
		\hline\hline
		$\mu>0$ & 1-2-3-4-5-6 & 1-2-3-4p-5-6p & 1-2-3-4-5-6  \\
		\hline
		$m(\sigma)=\sigma^{\gamma}$, $\gamma\geq 1$ & 1-2-3-4-5p & 
		1-2-3-4 & 
		1-2-3p-4p  \\
		\hline
		$m(\sigma)=\sigma^{\gamma}$, $\gamma\in(0,1)$ & 1-2-3-4 & 
		1-2-3-4 & 
		1-2-3p-4p  \\
		\hline
		$m(\sigma)\geq 0$ (Lip.\ cont.) & 1-2-3-4 & 1-2 & 1-2  \\
		\hline
	\end{tabular}
	\caption{state of the art on subproblems (1) through (6).  Numbers
	refer subproblems, a ``p'' means that in that case only partial
	(non-optimal) results have been obtained}
	\label{tab:overview}
\end{table}

Looking at Table~\ref{tab:overview}, one can guess that $p=1$ plays a
special role in the theory.  This is true also in the linear case.
Let us indeed consider equation
\begin{equation}
	au''(t)+\frac{b}{(1+t)^{p}}u'(t)+cAu(t)=0,
	\label{eq:wirth}
\end{equation}
where $a$, $b$, $c$ are positive parameters, and $p\geq 0$.  This
equation was investigated by T.\ Yamazaki~\cite{yamazaki-lin} and J.\
Wirth~\cite{wirth}.  They proved that (\ref{eq:wirth}) has both
parabolic and hyperbolic features, and which nature prevails depends
on $p$.  When $p<1$ the equation has parabolic behavior, in the sense
that all its solutions decay to 0 as $t\to +\infty$ as solutions of
the parabolic equation with $a=0$.  When $p>1$ the same equation has
hyperbolic behavior, meaning that every solution is asymptotic to a
suitable solution of the non-dissipative equation with $b=0$ (and in
particular all non-zero solutions do not decay to zero).  In the
critical case $p=1$ the nature of the problem depends on $b/a$, with
the parabolic behavior prevailing as soon as the ratio is large
enough.

Our results for Kirchhoff equations are consistent with the linear
theory.  Indeed we have always hyperbolic behavior when $p>1$, meaning
that non-zero global solutions (provided that they exist) cannot decay
to 0.  When $p\leq 1$ we were able to prove that the behavior is of
parabolic type in many cases.  In all such situations the critical
exponent $p=1$ falls in the parabolic regime, but this is simply due
to the fact that in our equation we have that $b=1$ and $a=\ep$ is
small enough, hence the ratio $b/a$ is always big enough.

For shortness's sake we don't include proofs in this note.
Nevertheless, we conclude this introduction by mentioning the useful
energies and the technical reasons why the problem becomes harder and
harder when the equation is degenerate, the dissipation is weak, and
the operator is non-coercive. In the parabolic case all estimates 
follow from the monotonicity of the classical energies
$$E_{k}(t):=|A^{k/2}u(t)|^{2},
\hspace{3em}
P(t):=\frac{|Au(t)|^{2}}{|A^{1/2}u(t)|^{2}}.$$

In the hyperbolic case, all known techniques for proving global
existence for (\ref{pbm:h-eq-gen}) require an a priori estimate such
as
\begin{equation}
	\ep\cdot\frac{\left|m'(|A^{1/2}\uep(t)|^{2})\right|}
	{m(|A^{1/2}\uep(t)|^{2})}
	\cdot|A\uep(t)|\cdot|\uep'(t)|\leq b(t).
	\label{est:basic}
\end{equation}

If $\mu>0$ and $b(t)$ is a positive constant, an a priori bound on
$A\uep(t)$ and $\uep'(t)$, together with the smallness of $\ep$, is
enough to establish (\ref{est:basic}).  When $b(t)\to 0$ as $t\to
+\infty$, the boundedness is no more enough, and we need some a priori
informations on the decay of $A\uep(t)$ and $\uep'(t)$.  This means
that global existence and decay estimates become intimately tied, and
they must be treated together.  The main energies involved in these
estimates are 
\begin{equation}
	E_{\ep,k}(t) := \ep\frac{|A^{k/2}\uep'(t)|^{2}}{\cep(t)}+
	|A^{(k+1)/2}\uep(t)|^{2},
	\hspace{2em}
	\Gep(t) := \frac{|\uep'(t)|^{2}}{{\cep^{2}(t)}},
	\label{defn:EG}
\end{equation}
where $\cep(t):=\m{\uep(t)}$.  They are both extensions of the first
energy of the parabolic case.  The use of $E_{\ep,k}(t)$ is quite
classical, and dates back to~\cite{debrito1,yamada}, while $\Gep(t)$ was
introduced by the authors in~\cite{gg:k-dissipative}.

The degenerate case is more complex.  Let us assume for example that
$m(\sigma)=0$ if and only if $\sigma=0$.  Then the decay of the
solution implies that the denominator in the left-hand side of
(\ref{est:basic}) tends to 0.  When $m(\sigma)=\sigma^{\gamma}$ with
$\gamma\in(0,1)$, then also the term with $m'$ diverges to $+\infty$
as the solution approaches 0.  This complicates proofs both in the
case of constant, and in the case of weak dissipation.

The basic idea to deal with degenerate nonlinear terms of the form
$m(\sigma)=\sigma^{\gamma}$ is to exploit that $\sigma
m'(\sigma)/m(\sigma)$ has a finite limit as $\sigma\to 0^{+}$. This 
reduces (\ref{est:basic}) to
\begin{equation}
	\ep\cdot\frac{|A\uep(t)|\cdot|\uep'(t)|}{|A^{1/2}\uep(t)|^{2}}
	\leq b(t).
	\label{est:basic+}
\end{equation}

This inequality has been approached using (\ref{defn:EG}) and the 
further energies
$$ \Pep  :=  \frac{\ep}{\cep} \frac{\au{\uep}^{2}\au{\uep'}^{2}-
\langle A\uep,\uep'\rangle^{2}}{\au{\uep}^{4}}+
\frac{|A\uep|^{2}}{\au{\uep}^{2}}, 
\hspace{3em}
\Qep  :=  \frac{|\uep'|^{2}}{\cep^{2}\auq{\uep}}.$$

These energies have been introduced by the first author
in~\cite{ghisi2} as a hyperbolic version of the second energy of the
parabolic case.

If $m(\sigma)$ is a general nonnegative (even Lipschitz continuous)
function, then it may happen that $\sigma m'(\sigma)/m(\sigma)$ is
unbounded in a neighborhood of $\sigma=0$.  In this case
(\ref{est:basic+}) does not imply (\ref{est:basic}), and each step
seems to require new ideas.  For this reason we are quite skeptic
about a future relevant progress in the last line of
Table~\ref{tab:overview}.

Concerning the coerciveness of the operator, it is well known that
small eigenvalues deteriorate the decay estimates on solutions, even
in the parabolic case, and we have seen that decay estimates are
fundamental also for existence issues.

This note is organized as follows.  Section~\ref{sec:parabolic} is
devoted to subproblem~(1), namely existence and decay estimates for
the parabolic problem.  Section~\ref{sec:local} is devoted to
subproblem~(2), namely local existence results for the hyperbolic
equation and local-in-time error estimates for the singular
perturbation problem.  In Section~\ref{sec:hyperbolic} we show that
(\ref{pbm:h-eq-gen}) has hyperbolic behavior whenever $p>1$.
Section~\ref{sec:ndg} is devoted to global existence and decay-error
estimates in the nondegenerate case.  Section~\ref{sec:dg} is devoted
to the degenerate case.  Finally, Section~\ref{sec:open} is a
collection of open problems.

\setcounter{equation}{0}
\section{The parabolic problem}\label{sec:parabolic}

The theory of parabolic equations of Kirchhoff type is quite well
established.  This equation appeared for the first time in the
pioneering paper~\cite{bernstein} by S.\ Bernstein.  He considered the
concrete equation in the interval $(0,1)$, with a nondegenerate
nonlinearity and constant dissipation, and he proved that for every
initial condition in the Sobolev space $H^{1}((0,1))$ the equation
admits a unique solution, which is actually analytic in the space
variable for every $t>0$ (the classical regularizing effect of
parabolic equations).  This result was afterwards extended by many
authors (see~\cite{bw,miletta}).

The more general version is probably stated in \cite{k-par}.  The basic
fact observed in \cite{k-par} is that any solution $u(t)$ of
(\ref{pbm:par}) can be written in the form
$$u(t)=v(\alpha(t)),$$
where $v(t)$ is the solution of the linear Cauchy problem with 
constant coefficients
\begin{equation}
	v'(t)+Av(t)=0,
	\quad\quad
	v(0)=u_{0},
	\label{eq:heat}
\end{equation}
and $\alpha:[0,+\infty)\to[0,+\infty)$ is the solution of the 
ordinary differential equation
$$b(t)\alpha'(t)=\m{v(\alpha(t))},
\quad\quad
\alpha(0)=0.$$

In other words, the solution of (\ref{pbm:par}) is always a time 
reparametrization of the solution of the heat-like equation 
(\ref{eq:heat}). At this point it is quite easy to prove the 
following existence result.

\begin{thm}[Global existence for the parabolic problem]
	Let $H$ be a Hilbert space, and let $A$ be a nonnegative
	self-adjoint (unbounded) operator on $H$ with dense domain.  Let
	$m:[0,+\infty)\to[0,+\infty)$ be a locally Lipschitz continuous
	function, and let $u_{0}\in\da$.
	
	Then problem (\ref{pbm:par}) has a unique global solution 
	$$u\in C^{1}\left([0,+\infty);H\right)\cap
	C^{0}\left([0,+\infty);\da\right).$$
	
	If in addition $\m{u_{0}}Au_{0}\neq 0$, hence $u'(0)\neq 0$, then
	the solution is non-stationary, and $u\in
	C^{1}\left((0,+\infty);D(A^{\alpha})\right)$ for every $\alpha\geq
	0$.
\end{thm}

Decay estimates for $u(t)$ can be deduced from decay estimates for
(\ref{eq:heat}) and the asymptotic behavior of the parametrization
$\alpha(t)$.  Concerning (\ref{eq:heat}), it is well known that the
asymptotic behavior of solutions depends on the coerciveness of the
operator $A$.  If $A$ is coercive with some constant $\nu>0$, then
solutions decay exponentially to 0, with a rate depending on $\nu$.
In this case we have indeed that 
$$|A^{1/2}u_{0}|^{2}\exp\left(-2
\frac{|Au_{0}|^{2}}{\auq{u_{0}}}\,t\right)\leq
|A^{1/2}v(t)|^{2}\leq|A^{1/2}u_{0}|^{2} \exp(-2\nu t).$$

If $A$ is non-coercive ($\nu\geq 0$), then decay rates are slower. We 
have indeed that
$$|A^{1/2}u_{0}|^{2}\exp\left(-2
\frac{|Au_{0}|^{2}}{\auq{u_{0}}}\,t\right)\leq\auq{v(t)}\leq
\frac{|u_{0}|^{2}}{2t},
\quad\quad
|Av(t)|^{2}\leq\frac{|u_{0}|^{2}}{2t^{2}}.$$

Note in particular that the estimates from below and from above for
$\auq{v(t)}$ involve different rates.  This range of rates cannot be
improved because, when the operator has a sequence of eigenvalues
converging to 0, any intermediate rate is realized by a suitable
solution.

Once we know the decay of $v(t)$, we can easily deduce the asymptotic
behavior of $\alpha(t)$, hence also the asymptotic behavior of 
$u(t)$. In Table~\ref{tab:decay-par} we sum up the decay estimates 
which can be obtained in this way, limiting ourselves for simplicity 
to  dissipation coefficients of the form $b(t)=(1+t)^{-p}$ with 
$p\geq 0$, and to nonlinear terms which are either non-degenerate or 
of the form $m(\sigma)=\sigma^{\gamma}$ with $\gamma>0$.

\begin{table}[htb]
	\renewcommand{\arraystretch}{1.3}
\begin{center}
	\begin{tabular}{|c||c|c|}
		\hline
		& $\nu>0$ & \allarga
		$\displaystyle{c_{1}e^{-\alpha_{1}(1+t)^{p+1}}\leq\auq{ u(t)}
		\leq c_{2}e^{-\alpha_{2}(1+t)^{p+1}}}$ \\
		& & \allarga
		$\displaystyle{c_{1}e^{-\alpha_{1}(1+t)^{p+1}}\leq
		|A u(t)|^{2} \leq c_{2}e^{-\alpha_{2}(1+t)^{p+1}}}$ \\
		& & \allarga $\displaystyle{c_{1}(1+t)^{2p}
		e^{-\alpha_{1}(1+t)^{p+1}} \leq| u'(t)|^{2}\leq
		c_{2}(1+t)^{2p} e^{-\alpha_{2}(1+t)^{p+1}}}$ \\
		\cline{2-3}
		\rotatebox{90}{\makebox[0pt]{$m(\sigma)\geq\mu>0$}} & $\nu\geq
		0$ & \allarga
		$\displaystyle{c_{1}e^{-\alpha_{1}(1+t)^{p+1}}\leq\auq{ u(t)}
		\leq \frac{c_{2}}{(1+t)^{p+1}}}$ \\
		& & \allarga
		$\displaystyle{|A u(t)|^{2} \leq \frac{c}{(1+t)^{2(p+1)}}}$ \\
		& & 	\allarga
		$\displaystyle{|u'(t)|^{2}\leq\frac{c}{(1+t)^{2}}}$ \\
		\hline
		\hline
		& $\nu>0$ & \allarga
		$\displaystyle{\frac{c_{1}}{(1+t)^{(p+1)/\gamma}}\leq\auq{ u(t)}
		\leq \frac{c_{2}}{(1+t)^{(p+1)/\gamma}}}$ \\
		&  & \allarga
		$\displaystyle{\frac{c_{1}}{(1+t)^{(p+1)/\gamma}}\leq
		|A u(t)|^{2} \leq \frac{c_{2}}{(1+t)^{(p+1)/\gamma}}}$ \\
		& & \allarga $\displaystyle{\frac{c_{1}}{(1+t)^{2+(p+1)/\gamma}}
		 \leq| u'(t)|^{2}\leq \frac{c_{1}}{(1+t)^{2+(p+1)/\gamma}}}$ \\
		 \cline{2-3}
		\rotatebox{90}{\makebox[0pt]{$m(\sigma)=\sigma^{\gamma}$}} &
		$\nu\geq 0$ & \allarga
		$\displaystyle{\frac{c_{1}}{(1+t)^{(p+1)/\gamma}}\leq\auq{ u(t)}
		\leq \frac{c_{2}}{(1+t)^{(p+1)/(\gamma+1)}}}$ \\
		& & \allarga $\displaystyle{|A u(t)|^{2} \leq
		\frac{c}{(1+t)^{(p+1)/\gamma}}}$ \\
		& & \allarga
		$\displaystyle{|u'(t)|^{2}\leq
		\frac{c}{(1+t)^{[2\gamma^{2}+(1-p)\gamma+p+1)/(\gamma^{2}+\gamma)}}}$
		\\
		\hline
	\end{tabular}
\caption{Decay estimates for the parabolic problem}
		\label{tab:decay-par}
		\end{center}
\end{table}

We stress that in all these cases solutions decay to zero, and the
decay rate becomes stronger and stronger as $p$ grows.  This contrasts
with the hyperbolic case, where solutions cannot decay when $p>1$ (see
Section~\ref{sec:hyperbolic}).

\setcounter{equation}{0}
\section{Local-in-time error estimates}\label{sec:local}

All the local existence results for the non-dissipative equation
(see~\cite[Theorem~2.1]{gg:survey-nd}) can be easily extended to the
dissipative case.  This provides a continuum of local existence
results, with the regularity requirements on the initial data
depending on the continuity modulus of $m$.  In this note we limit
ourselves to Lipschitz continuous nonlinear terms, or to the
non-Lipschitz case $m(\sigma)=\sigma^{\gamma}$ with $\gamma\in(0,1)$,
where the nondegeneracy assumption (\ref{hp:mdg}) makes the problem
just mildly non-Lipschitz.  In all these cases the equation is locally
well posed for initial data in Sobolev spaces.

In this section we focus on a property which is slightly stronger than
local existence, and which could be called \emph{almost global
existence}.  The first result is indeed that the life span of
$\uep(t)$ tends to $+\infty$ as $\ep\to 0^{+}$.

\begin{thm}[Hyperbolic problem: almost global 
	existence]\label{thm:ag-existence}
	Let $H$ be a Hil\-bert space, and let $A$ be a nonnegative
	self-adjoint (unbounded) operator on $H$ with dense domain.  Let
	$m:[0,+\infty)\to[0,+\infty)$ and $b:[0,+\infty)\to(0,+\infty)$ be
	two locally Lipschitz continuous functions.  Let us assume that
	$(u_{0},u_{1})\in\da\times\dau$ satisfy the non-degeneracy
	assumption (\ref{hp:mdg}), and let $T>0$.
	
	Then there exists $\ep_{0}>0$ such that for every
	$\ep\in(0,\ep_{0})$ problem (\ref{pbm:h-eq-gen}),
	(\ref{pbm:h-data}) has a unique solution
	$$\uep\in C^{2}\left([0,T];H\right)\cap
	C^{1}\left([0,T];\dau\right)\cap C^{0}\left([0,T];\da\right).$$
\end{thm}

Then we study the convergence of $\uep(t)$ to the solution $u(t)$ of
the limit problem.

\begin{thm}[Singular perturbation: local-in-time error estimates]
	\label{thm:lit-error}
	\hskip 0em plus 1em\mbox{} Let $H$, $A$, $m(\sigma)$, $b(t)$,
	$u_{0}$, $u_{1}$, $T$, $\ep_{0}$ be as in
	Theorem~\ref{thm:ag-existence}.  Let $u(t)$ be the solution of the
	corresponding parabolic problem (\ref{pbm:par}), and let $\rep(t)$
	and $\roep(t)$ be defined by (\ref{defn:rep}).
	
	Then we have the following conclusions.
	\begin{enumerate}
		\renewcommand{\labelenumi}{(\arabic{enumi})}
	\item  Without further assumptions on initial data, hence 
	$(u_{0},u_{1})\in\da\times\dau$, we have that
	$$|\roep(t)|^{2}+|A^{1/2}\roep(t)|^{2}+
	|A\roep(t)|^{2}+|\rep'(t)|^{2}\to 0
	\quad 
	\mbox{uniformly in }[0,T],$$
	$$\int_{0}^{T}|A^{1/2}\rep'(t)|^{2}\,dt\to 0.$$

	\item  If in addition we assume that 
	$(u_{0},u_{1})\in\dat\times\dau$, then there exists a 
	constant $C$ such that for every $\ep\in(0,\ep_{0})$ we have 
	that
	$$|\roep(t)|^{2}+|A^{1/2}\roep(t)|^{2}+ \ep|\rep'(t)|^{2} \leq
	C\ep^{2}
	\quad\quad 
	\forall t\in[0,T],$$
	$$\int_{0}^{T}|\rep'(t)|^{2}\,dt\leq C\ep^{2}.$$
	
	\item  If in addition we assume that 
	$(u_{0},u_{1})\in D(A^{2})\times\da$, then there exists a 
	constant $C$ such that for every $\ep\in(0,\ep_{0})$ we have 
	that
	$$|A\roep(t)|^{2}+|\rep'(t)|^{2}+\ep|A^{1/2}\rep'(t)|^{2}
	\leq C\ep^{2}\quad\quad\forall t\in[0,T],$$
	$$\int_{0}^{T}|A^{1/2}\rep'(t)|^{2}\,dt\leq C\ep^{2}.$$
	\end{enumerate}
\end{thm}

We point out that the remainder $\rep(t)$ is well suited for estimates
involving derivatives, because it doesn't feel the effects of the
boundary layer due to the loss of one initial condition.  On the
contrary, the remainder $\roep(t)$ is better suited for estimates
without derivatives.  This is because, for example, $A\roep(0)$ is
defined whenever $u_{0}\in D(A)$, while $A\rep(0)$ requires $u_{0}\in
D(A^{2})$ (see definition~(\ref{pbm:theta-data}) of $w_{0}$).

Both the existence and the convergence result are local-in-time,
namely constants, error estimates, and the smallness of $\ep$ do
depend on the interval $[0,T]$ chosen at the beginning.  On the other
hand, the assumptions required on $b(t)$ and $m(\sigma)$ are quite
weak.  The dichotomy between hyperbolic and parabolic behavior
mentioned in the introduction appears only as $t\to +\infty$, hence it
plays no role on a fixed time interval.  In particular we don't need
to assume that $p\leq 1$ in the case where $b(t)=(1+t)^{-p}$.

From Theorem~\ref{thm:lit-error} it is clear that convergence rates
for the singular perturbation problem depend on the regularity of
initial data. This situation is consistent with the linear case. 
Indeed in~\cite{gg:l-cattaneo} we considered the linear equation
$$\ep\uep''(t)+\uep'(t)+A\uep(t)=0$$
and the corresponding limit parabolic problem, and we proved similar
results.  We also proved that an error estimate such as
$|A^{1/2}\roep(t)|^{2}\leq C\ep^{2}$ is possible only when
$(u_{0},u_{1})\in\dat\times\dau$.

We just remark that uniform convergence, without any rate, requires
initial data in spaces such as $\da\times\dau$, hence with ``gap 1/2''
between the regularity of $\uep'$ and $\uep$, a typical feature of
hyperbolic problems.  On the contrary, if we want some convergence
rate, we have to work in spaces such as $\dat\times\dau$ or
$D(A^{2})\times\da$, hence with ``gap 1'', a typical feature of
parabolic problems.  In our opinion, this gives further evidence that
the parabolic nature dominates in the limit.

A formal proof of Theorem~\ref{thm:ag-existence} and
Theorem~\ref{thm:lit-error}, as they are stated, has never been put
into writing.  Error estimates have been considered in at least three
papers, always with constant dissipation.  B.\ F.\ Esham, and R.\ J.\
Weinacht~\cite{ew} proved error estimates in the nondegenerate case
with initial data in $\dat\times\da$.  The second
author~\cite{k-cattaneo} considered the degenerate case, proving
uniform convergence for data in $\da\times\dau$, and error estimates
for more regular data.  Finally, the authors~\cite{gg:k-PS} proved
error estimates in the degenerate case with the optimal requirement
that initial data are in $\dat\times\dau$
(see~\cite[Proposition~A.1]{gg:k-PS}).  It should be quite standard to
extend those proofs to the case of weak dissipation, just because on a
fixed time interval the function $b(t)$ is always strictly positive.

\setcounter{equation}{0}
\section{The hyperbolic regime}\label{sec:hyperbolic}

In this section we show that when the dissipation is too weak, namely
$b(t)\to 0$ too fast, then equation (\ref{pbm:h-eq-gen}) behaves in a
hyperbolic way, in the sense that its non-zero global solutions
(provided that they exist) do not decay to 0 as $t\to+\infty$.  Of
course this doesn't prevent such solutions from existing (which
remains an open problem), but it shows that the problem cannot be
approached using the standard methods based on estimates such as
(\ref{est:basic}) or (\ref{est:basic+}).  Since solutions of the
parabolic problem always decay to 0, this shows also that no
decay-error estimate can be true in this case.  Note that condition
(\ref{hp:int-b}) is equivalent to $p>1$ when $b(t)=(1+t)^{-p}$.

\begin{thm}[Hyperbolic regime]\label{thm:p>1}
	Let $H$ be a Hilbert space, and let $A$ be a nonnegative
	self-adjoint (unbounded) operator on $H$ with dense domain.  
	
	Let
	$m:[0,+\infty)\to[0,+\infty)$ be a continuous function.  Let
	$b:[0,+\infty)\to(0,+\infty)$ be a continuous function such that
	\begin{equation}
		\int_{0}^{+\infty}b(s)\,ds<+\infty.
		\label{hp:int-b}
	\end{equation}
	
	Let $(u_{0},u_{1})\in\da\times\dau$ be such that
	\begin{equation}
		|u_{1}|^{2}+\int_{0}^{|A^{1/2}u_{0}|^{2}}m(\sigma)\,d\sigma>0.
		\label{hp:p>1}
	\end{equation}
	
	Let us assume that for some $\ep>0$ problem (\ref{pbm:h-eq-gen}),
	(\ref{pbm:h-data}) has a global solution 
	\begin{equation}
		\uep\in C^{2}([0,+\infty);H)\cap C^{1}([0,+\infty);\dau) \cap
		C^{0}([0,+\infty);\da).
		\label{th:uep-reg}
	\end{equation}	
	
	Then we have that
	\begin{equation}
		\liminf_{t\to +\infty}\left(
		|\uep'(t)|^{2}+\auq{\uep(t)}\right)>0.
		\label{th:p>1}
	\end{equation}
\end{thm}

The proof of this result is very simple, and relies on the usual 
Hamiltonian
$$\Hep(t):=\ep|\uep'(t)|^{2}+
\int_{0}^{|A^{1/2}\uep(t)|^{2}}m(\sigma)\,d\sigma.$$

Assumption (\ref{hp:p>1}) is equivalent to say that $\Hep(0)>0$.
Moreover we have that
$$\Hep'(t)=-2b(t)|\uep'(t)|^{2}\geq -\frac{2}{\ep}b(t)\Hep(t)
\quad\quad\forall t\geq 0,$$
hence
$$\Hep(t)\geq \Hep(0)\exp\left(-\frac{2}{\ep}\int_{0}^{t}b(s)\,ds\right)
\quad\quad\forall t\geq 0.$$

For a fixed $\ep>0$, the right-hand side is greater than a positive
constant independent on $t$ because of (\ref{hp:int-b}) and the fact
that $\Hep(0)>0$.  This implies (\ref{th:p>1}).

\setcounter{equation}{0}
\section{The nondegenerate case}\label{sec:ndg}

In this section we focus on the hyperbolic equation
(\ref{pbm:h-eq-gen}) under the non-degeneracy assumption $\mu>0$.

The case with constant dissipation was considered independently by E.\
H.\ de Brito~\cite{debrito1} and by Y.\ Yamada~\cite{yamada}.  They
proved existence of a global solution provided that $\ep$ is small
enough.  Decay estimates for these solutions were proved by Y.\
Yamada~\cite{yamada} in the non-coercive case, and by E.\ H.\ de
Brito~\cite{debrito2} and by M.\ Hosoya and Y.\ Yamada~\cite{hy} in
the coercive case.  All these estimates were afterwards reobtained as
a particular case of the theory developed in~\cite{gg:k-decay}.

More recently, H.\ Hashimoto and T.\ Yamazaki~\cite{yamazaki} proved
that for initial data $(u_{0},u_{1})\in D(A^{3/2})\times D(A)$ one has
that
$$|\roep(t)|^{2}+(1+t)|A^{1/2}\roep(t)|^{2}+
\ep(1+t)^{2}|\rep'(t)|^{2}\leq C\ep^{2} 
\quad\quad
\forall t\geq 0,$$
where of course $C$ doesn't depend on $t$ and $\ep$.  When
$(u_{0},u_{1})\in D(A^{2})\times D(A)$, the coefficient $\ep$ in the
left-hand side may be dropped, thus providing a better convergence
rate on $\rep'(t)$.  This is a first example of \emph{decay-error
estimate}.

The weakly dissipative case was considered only in last years.  Apart
from a result obtained in a special situation by M.\ Nakao and J.\
Bae~\cite{nakao}, the problem in its full generality was solved by T.\
Yamazaki~\cite{yamazaki-wd} in the subcritical case $p<1$ with some
technical requirements on initial data, and then by the
authors~\cite{gg:w-ndg} (see also~\cite{yamazaki-cwd}) in the general
case $p\leq 1$ with minimal requirements on initial data.

The results are the following.

\begin{thm}[Hyperbolic problem: global existence]\label{thm:h-existence} 
	Let $H$ be a Hilbert space, and let $A$ be a nonnegative
	self-adjoint (unbounded) operator on $H$ with dense domain.  Let
	$\mu>0$, and let $m:[0,+\infty)\to[\mu,+\infty)$ be a locally
	Lipschitz continuous function.  Let $b(t):=(1+t)^{-p}$ with
	$p\in[0,1]$, and let $(u_{0},u_{1})\in\da\times\dau$.
			
	Then there exists $\ep_{0}>0$ such that for every
	$\ep\in(0,\ep_{0})$ problem (\ref{pbm:h-eq-gen}),
	(\ref{pbm:h-data}) has a unique global solution $\uep$ satisfying
	(\ref{th:uep-reg}).
\end{thm}

\begin{thm}[Hyperbolic problem: decay estimates]\label{thm:h-decay} 
	Under the same assumptions of Theorem~\ref{thm:h-existence} there
	exists a constant $C$ such that for every $\ep\in(0,\ep_{0})$ we 
	have that
	$$|\uep(t)|^{2}+(1+t)^{p+1}\auq{\uep(t)}+(1+t)^{2}|\uep'(t)|^{2}\leq C
		\quad\quad\forall t\geq 0,$$
	$$\ep|A^{1/2}\uep'(t)|^{2}+ |A\uep(t)|^{2}\leq 
		\frac{C}{(1+t)^{2(p+1)}}
		\quad\quad\forall t\geq 0,$$
	$$\int_{0}^{+\infty}(1+t)^{p}\left(
			|\uep'(t)|^{2}+\auq{\uep(t)}\right)\,dt\leq C,$$
	$$\int_{0}^{+\infty}(1+t)^{2p+1}\left(
			|A^{1/2}\uep'(t)|^{2}+|A\uep(t)|^{2}\right)\,dt\leq C.$$
\end{thm}

\begin{thm}[Singular perturbation: decay-error estimates]\label{thm:decay-error}
	\hskip 0em plus 1em\mbox{}
	Let $H$, $A$, $\mu$, $m(\sigma)$, $b(t)$, $p$, $u_{0}$, $u_{1}$,
	$\ep_{0}$ be as in Theorem~\ref{thm:h-existence}.  Let $u(t)$ be
	the solution of the corresponding parabolic problem
	(\ref{pbm:par}), and let $\rep(t)$ and $\roep(t)$ be defined by
	(\ref{defn:rep}).
	
	Then we have the following conclusions.
	\begin{enumerate}
		\renewcommand{\labelenumi}{(\arabic{enumi})}
		\item  Without further assumptions on initial data, namely 
		$(u_{0},u_{1})\in\da\times\dau$, we have that
		$$|\roep(t)|^{2}+(1+t)^{p+1}|A^{1/2}\roep(t)|^{2}+
		(1+t)^{2(p+1)}|A\roep(t)|^{2}+(1+t)^{2}|\rep'(t)|^{2}\to 0$$
		uniformly in $[0,+\infty)$, and
		$$\int_{0}^{+\infty}(1+t)^{p}\left(
		|\rep'(t)|^{2}+|A^{1/2}\roep(t)|^{2}\right)dt\to 0,$$
		$$\int_{0}^{+\infty} (1+t)^{2p+1}\left(
		|A^{1/2}\rep'(t)|^{2}+|A\roep(t)|^{2}\right)dt\to 0.$$
	
		\item  If in addition we assume that 
		$(u_{0},u_{1})\in\dat\times\dau$, then there exists a 
		constant $C$ such that for every $\ep\in(0,\ep_{0})$ we have 
		that
		$$|\roep(t)|^{2}+(1+t)^{p+1}|A^{1/2}\roep(t)|^{2}+
		\ep(1+t)^{p+1}|\rep'(t)|^{2}
		\leq C\ep^{2}
		\quad\quad
		\forall t\geq 0,$$
		$$\int_{0}^{+\infty}(1+t)^{p}\left(
		|\rep'(t)|^{2}+|A^{1/2}\roep(t)|^{2}\right)dt\leq C\ep^{2}.$$
		
		\item  If in addition we assume that 
		$(u_{0},u_{1})\in D(A^{2})\times\da$, then there exists a 
		constant $C$ such that for every $\ep\in(0,\ep_{0})$ we have 
		that
		$$(1+t)^{2(p+1)}|A\roep(t)|^{2}+(1+t)^{2}|\rep'(t)|^{2}
		\leq C\ep^{2}
		\quad\quad
		\forall t\geq 0,$$
		$$\int_{0}^{+\infty} (1+t)^{2p+1}\left(
		|A^{1/2}\rep'(t)|^{2}+|A\roep(t)|^{2}\right)dt\leq C\ep^{2}.$$
	\end{enumerate}
	
\end{thm}

We point out that in the previous three theorems the operator $A$ is 
never assumed to be coercive.

The global-in-time convergence rates (with respect to $\ep$) 
appearing in Theorem~\ref{thm:decay-error} are optimal because they 
coincide with the local-in-time convergence rates of 
Theorem~\ref{thm:lit-error}, which in turn are the same of the linear 
case.

The decay rates (with respect to time) appearing in 
Theorem~\ref{thm:h-decay} and Theorem~\ref{thm:decay-error} are 
optimal for non-coercive operators. In this case indeed they coincide 
with the decay rates of solutions of the corresponding parabolic 
equation, as shown in Table~\ref{tab:decay-par}.

In the coercive case these decay rates are not optimal.  In the case
$p=0$ indeed we know that solutions exponentially decay to zero as
solutions of the parabolic problem
(see~\cite{debrito2,hy,gg:k-decay}).  We strongly suspect that the
same is true also for every $p\in[0,1]$, namely that solutions decay
as shown in the first three rows of Table~\ref{tab:decay-par}.  Of
course also the decay rates in Theorem~\ref{thm:decay-error} should be
changed accordingly.  We give no precise statement or reference
because this part of the theory has never been put into writing.

\setcounter{equation}{0}
\section{The degenerate case}\label{sec:dg}

Several papers have been devoted to global existence and decay
estimates for equation (\ref{pbm:h-eq-gen}) in the degenerate case
$\mu\geq 0$.  Let us begin with constant dissipation.  In this case
global existence results (provided that the problem is mildly
degenerate and $\ep$ is small enough) were proved by K.\ Nishihara and
Y.\ Yamada~\cite{ny} in the case where $m(\sigma)=\sigma^{\gamma}$
(with $\gamma\geq 1$), by the authors~\cite{gg:k-dissipative} in the
case where $m(\sigma)\geq 0$ is any Lipschitz continuous function, and
by the first author~\cite{ghisi1,ghisi2} in the non-Lipschitz case
where $m(\sigma)=\sigma^{\gamma}$ with $\gamma\in(0,1)$.

Decay estimates have long been studied for equations with constant
dissipation.  In the case $m(\sigma)=\sigma^{\gamma}$ with $\gamma\geq
1$, the first decay estimates were obtained by K.\ Nishihara and Y.\
Yamada~\cite{ny} in the coercive case, and by K.\ Ono~\cite{ono-mm} in
the non-coercive case.  The case $m(\sigma)=\sigma^{\gamma}$ with
$\gamma\in(0,1)$ was considered in~\cite{ghisi1}.  In the special case
$m(\sigma)=\sigma$, T.\ Mizumachi~\cite{mizu-ade,mizu-nc} and K.\
Ono~\cite{ono-kyushu,ono-aa} proved better decay estimates, namely
estimates with decay rates which are faster than those obtained by
putting $\gamma=1$ in the previous ones.  This in particular showed
that the previous results were not optimal.

A complete answer was given by the authors in~\cite{gg:k-decay}, where
the case of a general nonlinearity $m(\sigma)\geq 0$ is considered.
The decay rates obtained in~\cite{gg:k-decay} coincide with the decay
rates of solutions of the parabolic problem.

Let us consider now the equation with weak dissipation, focussing on
the model case
\begin{equation}
	\ep\uep''(t)+\frac{1}{(1+t)^{p}}\uep'(t)+
	|A^{1/2}\uep(t)|^{2\gamma}A\uep(t)=0,
	\label{eq:model-w-dg}
\end{equation}
of course with the mild non-degeneracy assumption (\ref{hp:mdg}).  The
only previous result we are aware of was obtained by K.\
Ono~\cite{ono-wd}.  In the special case $\gamma=1$ he proved that a
global solution exists provided that $\ep$ is small and 
$p\in[0,1/3]$. The reason of the slow progress in this field is 
hardly surprising. In the weakly dissipative case existence and decay 
estimates must be proved in the same time. The better are the decay 
estimates, the stronger is the existence result.

Ten years ago decay estimates for degenerate equations were far from 
being optimal, but for the special case $\gamma=1$. 
In~\cite{gg:k-decay} a new method for obtaining optimal decay 
estimates was introduced. This allowed a substantial progress on 
equation (\ref{eq:model-w-dg}).

Let us begin with our existence and decay results proved in~\cite{gg:w-dg}.
The first one concerns the coercive case.

\begin{thm}[Coercive case: global existence and decay estimates]
	\label{thm:c-existence}
	Let $H$ be a Hilbert space, and let $A$ be a nonnegative
	self-adjoint (unbounded) operator with dense domain.  Let us assume
	that $A$ is coercive ($\nu>0$). Let $\gamma>0$, and let
	$p\in[0,1]$.  Let us assume that $(u_{0},u_{1})\in\da\times\dau$
	satisfy (\ref{hp:mdg}).
	
	Then there exists $\ep_{0}>0$ such that for every
	$\ep\in(0,\ep_{0})$ problem (\ref{eq:model-w-dg}),
	(\ref{pbm:h-data}) has a unique global solution satisfying 
	(\ref{th:uep-reg}).
	
	Moreover there exist positive constants $C_{1}$ and $C_{2}$ such that
	$$\frac{C_{1}}{(1+t)^{(p+1)/\gamma}}\leq
		\auq{\uep(t)}\leq
		\frac{C_{2}}{(1+t)^{(p+1)/\gamma}}
		\quad\quad\forall t\geq 0,$$
	$$\frac{C_{1}}{(1+t)^{(p+1)/\gamma}}\leq
		|A\uep(t)|^{2}\leq
		\frac{C_{2}}{(1+t)^{(p+1)/\gamma}}
		\quad\quad\forall t\geq 0,$$
	$$|\uep'(t)|^{2}\leq
		\frac{C_{2}}{(1+t)^{2+(p+1)/\gamma}}
		\quad\quad\forall t\geq 0.$$
\end{thm}

We point out that Theorem~\ref{thm:c-existence} is optimal both in the
sense that all $p\in[0,1]$ are considered, and in the sense that
solutions decay as in the parabolic case (see
Table~\ref{tab:decay-par}).

In the non-coercive case we have the following result.

\begin{thm}[Non-coercive case: global existence and decay estimates]
	\label{thm:nc-existence}
	\hskip 0em plus 1em\mbox{}
	Let $H$ be a Hilbert space, and let $A$ be a nonnegative
	self-adjoint (unbounded) operator with dense domain.  Let
	$\gamma\geq 1$, and let
	\begin{equation}
		0\leq p\leq \frac{\gamma^{2}+1}{\gamma^{2}+2\gamma-1}.  
		\label{hp:pnc}
	\end{equation}
	
	Let us assume
	that $(u_{0},u_{1})\in\da\times\dau$ satisfy (\ref{hp:mdg}).
	
	Then there exists $\ep_{0}>0$ such that for every
	$\ep\in(0,\ep_{0})$ problem (\ref{eq:model-w-dg}),
	(\ref{pbm:h-data}) has a unique global solution satisfying
	(\ref{th:uep-reg}).
	
	Moreover there exist constants $C_{1}$ and $C_{2}$ such that
	$$\frac{C_{1}}{(1+t)^{(p+1)/\gamma}}\leq
		\auq{\uep(t)}\leq
		\frac{C_{2}}{(1+t)^{(p+1)/(\gamma+1)}}
		\quad\quad\forall t\geq 0,$$
	$$|A\uep(t)|^{2}\leq
		\frac{C_{2}}{(1+t)^{(p+1)/\gamma}}
		\quad\quad\forall t\geq 0,$$
	$$|\uep'(t)|^{2}\leq
		\frac{C_{2}}{(1+t)^{[2\gamma^{2}+(1-p)\gamma+p+1]/
		(\gamma^{2}+\gamma)}}
		\quad\quad\forall t\geq 0.$$
\end{thm}

Theorem~\ref{thm:nc-existence} doesn't represent a final answer in the
non-coercive case.  Let indeed $p_{\gamma}$ denote the right-hand side
of (\ref{hp:pnc}).  It is easy to see that $p_{\gamma}\leq 1$ for
every $\gamma\geq 1$, with equality only when $\gamma=1$, and
asymptotically as $\gamma\to +\infty$.  Since we have hyperbolic
behavior when $p> 1$ (see Section~\ref{sec:hyperbolic}), and parabolic
behavior for $p\in[0,p_{\gamma}]$, this means that there is a
non-man's land between $p_{\gamma}$ and 1 where things are not clear
yet.

The only case where this region is empty is when $\gamma=1$.  In this
case all exponents $p\in[0,1]$ fall in the parabolic regime, and this
improves the result obtained in \cite{ono-wd} ($p\in[0,1/3]$) also in
the case $m(\sigma)=\sigma$.

\begin{figure}[htbp]
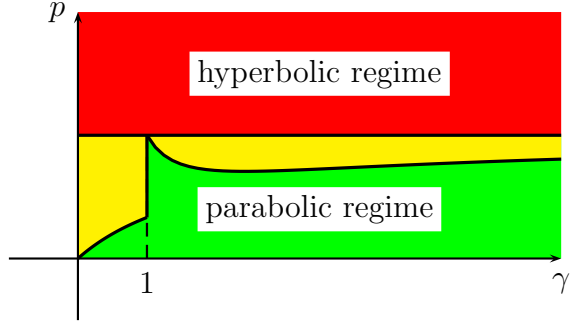

	\centering
	\psset{xunit=5ex,yunit=9ex}
	\pspicture(-1,-0.5)(7,2)
	\SpecialCoor
	\psframe*[linecolor=red](0,1)(7,2)
	\psframe*[linecolor=green](0,0)(7,1)
	\pscustom[fillcolor=yellow,fillstyle=solid,linestyle=none]{
	\psplot{0}{1}{x x 2 add div}
	\psline(!1 1 3 div)(1,1)
	\psplot{1}{7}{x 2 exp 1 add x 2 exp 2 x mul 1 sub add div}
	\psline(!7 50 62 div)(7,1)
	\psline(7,1)(0,1)
	\psline(0,1)(0,0)
	}
	\psline{->}(-1,0)(7,0)
	\psline{->}(0,-0.5)(0,2)
	\psset{linewidth=1.5\pslinewidth}
	\psplot{0}{1}{x x 2 add div}
	\psline(!1 1 3 div)(1,1)
	\psplot{1}{7}{x 2 exp 1 add x 2 exp 2 x mul 1 sub add div}
	\psline(7,1)(0,1)
	\psline[linestyle=dashed,linewidth=0.7\pslinewidth](1,0)(1,1)
	\rput*(3.5,1.5){hyperbolic regime}
	\rput*(3.5,0.4){parabolic regime}
	\rput(7,-0.2){$\gamma$}
	\rput(-0.3,2){$p$}
	\rput(1,-0.2){$1$}
	\endpspicture
	\caption{parabolic and hyperbolic regimes in the degenerate non-coercive case}
	\label{fig:no-man-land}
\end{figure}

We stated Theorem~\ref{thm:nc-existence} assuming $\gamma\geq 1$.  In
the case $\gamma\in (0,1)$ we have a weaker result, namely global
existence for $p\in[0,\gamma/(\gamma+2)]$
(see~\cite[Remark~2.6]{gg:w-dg}). Figure~\ref{fig:no-man-land} 
represents hyperbolic and parabolic regimes, and the no-man's land in 
between.

The singular perturbation problem is still quite open in the
degenerate case.  We have indeed only the following partial result for
the constant dissipation case (see~\cite{gg:k-PS}).

\begin{thm}[Constant dissipation: global-in-time error estimates]
	\label{thm:dg-error}
	Let $H$ be a Hilbert space, and let $A$ be a nonnegative
	self-adjoint (unbounded) operator with dense domain.  Let $\uep(t)$
	be the solution of equation (\ref{eq:model-w-dg}) with $\gamma>0$,
	$p=0$, and initial data $(u_{0},u_{1})\in\da\times\dau$ satisfying
	(\ref{hp:mdg}).  Let $u(t)$ be the solution of the corresponding
	parabolic problem, and let $\rep(t)$ and $\roep(t)$ be defined by
	(\ref{defn:rep}).
	
	Then we have the following conclusions.
	\begin{enumerate}
		\renewcommand{\labelenumi}{(\arabic{enumi})} 
		\item Without
		further assumptions on initial data, namely
		$(u_{0},u_{1})\in\da\times\dau$, we have that
		$$|\roep(t)|^{2}+|A^{1/2}\roep(t)|^{2}+
		|A\roep(t)|^{2}+|\rep'(t)|^{2}\to 0$$
		uniformly in $[0,+\infty)$, and 
		$$\int_{0}^{+\infty}\left(
		|\rep'(t)|^{2}+|A^{1/2}\rep'(t)|^{2}\right)dt\to 0.$$
	
		\item  If in addition we assume that $\gamma\geq 1$ and
		$(u_{0},u_{1})\in\dat\times\dau$, then there exists a 
		constant $C$ such that for every $\ep\in(0,\ep_{0})$ we have 
		that
		$$|\roep(t)|^{2}+\ep|A^{1/2}\roep(t)|^{2}\leq C\ep^{2}
		\quad\quad
		\forall t\geq 0,$$
		$$\int_{0}^{+\infty}
		|\rep'(t)|^{2}\,dt\leq C\ep.$$
		
		\item  If in addition we assume that $\gamma\geq 1$ and 
		$(u_{0},u_{1})\in\dat\times\da$, then there exists a 
		constant $C$ such that for every $\ep\in(0,\ep_{0})$ we have 
		that
		$$|\roep(t)|^{2}+\ep^{2/3}|A^{1/2}\roep(t)|^{2}+
		\ep^{4/3}|A\roep(t)|^{2}+\ep^{4/3}|\rep'(t)|^{2}
		\leq C\ep^{2}
		\quad\quad
		\forall t\geq 0.$$
	\end{enumerate}
\end{thm}

Theorem~\ref{thm:dg-error} is far from being optimal.  First of all
most of the convergence rates in the second and third statement are
weaker than the corresponding rates in Theorem~\ref{thm:lit-error}.
Moreover, all statements present just error estimates, and not
decay-error estimates as in Theorem~\ref{thm:decay-error}.  It is
possible to add some decays with some extra work, but in any case they
are so far from those appearing in Theorem~\ref{thm:c-existence} and
Theorem~\ref{thm:nc-existence} that we decided not to include them.
Last but not least, Theorem~\ref{thm:dg-error} is limited to equations
with constant dissipation.

\setcounter{equation}{0}
\section{Open problems}\label{sec:open}

The main open problem in the theory of Kirchhoff
equations is existence of global solutions. We have seen that in the 
dissipative case an affirmative answer can be given provided that 
$\ep$ is small enough. So the first question is whether this 
condition is necessary or not.

\begin{open}
	Let us consider equation (\ref{pbm:h-eq-gen}) with
	$m:[0,+\infty)\to[1,+\infty)$ of class $C^{\infty}$, and constant
	dissipation $b(t)\equiv 1$.  Let us assume that $(u_{0},u_{1})\in
	D(A^{\infty})\times D(A^{\infty})$, where $D(A^{\infty})$ is the
	intersection of all spaces $D(A^{\alpha})$ with $\alpha\geq 0$.
	
	Does the Cauchy problem (\ref{pbm:h-eq-gen}), (\ref{pbm:h-data}) 
	admit a global solution for every $\ep>0$?
\end{open}

We stated the question with generous assumptions both on the
nonlinearity (smoothness and strict hyperbolicity), and on initial
data (regularity). In any case there are no counterexamples, even 
with less regular terms and data, or with $b(t)\equiv 0$.

Even assuming the smallness of $\ep$, one may ask if a global 
solution exists under assumptions weaker than those required in the 
previous sections. This leads to the following question.

\begin{open}
	Let us consider the Cauchy problem (\ref{pbm:h-eq-gen}), 
	(\ref{pbm:h-data}) in each of the following situations.
	\begin{itemize}
		\item  In the hyperbolic regime where $b(t)=(1+t)^{-p}$ with 
		$p>1$.
	
		\item In the case where assumption $b(t)=(1+t)^{-p}$ with
		$p\leq 1$ is replaced by the weaker condition
		(\ref{hp:int-b}).
	
		\item  In the really degenerate case $\m{u_{0}}=0$.
	\end{itemize}
	
	Is it possible to prove global existence provided that $\ep$ is 
	small enough?
\end{open}

A third question related to global existence issues concerns the
regularity of initial data.  All existence results stated in the
previous sections assume that $(u_{0},u_{1})\in\da\times\dau$.  On the
other hand, the classical local existence results for the
non-dissipative equation require the weaker assumption $(u_{0},u_{1})\in
D(A^{3/4})\times D(A^{1/4})$.  Therefore a natural question is whether
the global existence results for dissipative equations can be extended
to this weaker class of data.  

In the constant dissipation case, it is not difficult to give an
affirmative answer when $\mu>0$ or when $m(\sigma)=\sigma^{\gamma}$
with $\gamma\geq 2$.  On the contrary, the proof given in
\cite{gg:k-dissipative} for a general locally Lipschitz continuous
non-linearity $m(\sigma)\geq 0$ seems to require in an essential way
that $(u_{0},u_{1})\in\da\times\dau$.  So the problem is the
following.

\begin{open}
	
	Let us consider equation (\ref{pbm:h-eq-gen}) with constant
	dissipation $b(t)\equiv 1$, and with any locally Lipschitz
	continuous nonlinearity $m(\sigma)\geq 0$.  Let us assume that
	$(u_{0},u_{1})\in D(A^{3/4})\times D(A^{1/4})$ satisfy the
	non-degeneracy condition (\ref{hp:mdg}).
	
	Does problem (\ref{pbm:h-eq-gen}), (\ref{pbm:h-data}) admit a global
	solution for every small enough $\ep$?
\end{open}

The last open question concerning existence is how to fill the
no-man's zone left by Theorem~\ref{thm:nc-existence} and described in
Figure~\ref{fig:no-man-land}.

\begin{open}
	Let $H$, $A$, $u_{0}$, $u_{1}$ be as in
	Theorem~\ref{thm:nc-existence}.  Let us assume that either
	$\gamma\in(0,1)$ and $p\in(\gamma/(\gamma+2),1]$, or that
	$\gamma>1$ and $p\in(p_{\gamma},1]$, where $p_{\gamma}$ is the
	right-hand side of (\ref{hp:pnc}).
	
	Does problem (\ref{eq:model-w-dg}), (\ref{pbm:h-data}) admit a 
	global solution whenever $\ep$ is small enough?
\end{open}

All previous examples suggest that the answer should be affirmative, 
but a proof seems to require some new ideas.

The singular perturbation problem is arguably the new frontier in this
research field.  This problem has been quite well understood only in the
nondegenerate case, in which case, however, the decay rates are
optimal only for non-coercive operators. A first open question is 
therefore the following.

\begin{open}
	Let the assumptions of Theorem~\ref{thm:decay-error} be 
	satisfied. Let us assume also that the operator $A$ is coercive 
	($\nu>0$).
	
	Prove the same conclusions of Theorem~\ref{thm:decay-error} with 
	all polynomial decay rates such as $(1+t)^{\beta}$ replaced by 
	exponential decay rates of the form $\exp(\alpha(1+t)^{p+1})$, 
	where $\alpha$ is a suitable constant.
\end{open}

The singular perturbation problem is quite open in the degenerate
case.  One should try to extend Theorem~\ref{thm:dg-error} in order to
allow weak dissipations, and involve better decay and convergence
rates.  An example of open question is the following.

\begin{open}
	Let $H$, $A$, $\gamma$, $p$, $u_{0}$, $u_{1}$, $\ep_{0}$ be as in
	Theorem~\ref{thm:c-existence}.  Let $u(t)$ be the solution of the
	corresponding parabolic problem, and let $\rep(t)$ and $\roep(t)$
	be defined by (\ref{defn:rep}).
	
	Under the appropriate conditions on initial data, prove that
	there exists a constant $C$ such that for every
	$\ep\in(0,\ep_{0})$ we have that
	$$(1+t)^{(p+1)/\gamma}|A^{1/2}\roep(t)|^{2}\leq C\ep^{2}
	\quad\quad
	\forall t\geq 0,$$
	$$(1+t)^{2+(p+1)/\gamma}|\rep'(t)|^{2}\leq C\ep^{2}
	\quad\quad
	\forall t\geq 0.$$
\end{open}

In this estimates we require on $\rep(t)$ and $\roep(t)$ the same decay
rates (as $t\to +\infty$) of $\uep(t)$ and $u(t)$ separately, and we
require the same convergence rates (as $\ep\to 0^{+}$) of the
local-in-time error estimates.  We actually suspect that in the
degenerate case the remainders $\rep(t)$ and $\roep(t)$ decay faster
than $\uep(t)$ and $u(t)$.

\subsubsection*{\centering Acknowledgments}

This note is an extended version of the talk presented by the first
author in the section ``Dispersive Equations'' of the 7th ISAAC
conference (London 2009).  We would like to thank once again the
organizers of that section, Prof.\ F.\ Hirosawa and Prof.\ M.\
Reissig, for their kind invitation. We would like to thank also 
Prof.\ T.\ Yamazaki for sending us preliminary versions of 
references~\cite{yamazaki-wd,yamazaki-cwd}, and for pointing out 
reference~\cite{wirth}.

\label{NumeroPagine}

\end{document}